\documentclass[a4paper, reqno, 12pt]{amsart}
\usepackage[utf8]{inputenc} 
\usepackage[T1, T2A]{fontenc}
\usepackage[english]{babel}
\usepackage{amsmath, amssymb, amsfonts, amsthm, amscd, mathrsfs}
\usepackage{enumitem}
\usepackage[hidelinks]{hyperref}
\usepackage[usenames]{color}
\setlist[enumerate]{label = (\alph*), ref=(\text{\alph*)}}
\setlist[itemize]{nolistsep}

\usepackage{epsfig}
\usepackage{graphicx}
\usepackage[cmtip,arrow]{xy}
\usepackage{bm}
\usepackage{tikz}
\usetikzlibrary{calc}
\usetikzlibrary[shapes.geometric]
\usepackage{multirow}
\usepackage{comment}

\renewcommand{\phi}{\varphi}
\renewcommand{\ge}{\geqslant}
\renewcommand{\le}{\leqslant}

\newcommand{\KK}{\mathbb{K}}
\renewcommand{\AA}{\mathbb{A}}

\newcommand{\ZZ}{\mathbb{Z}}
\newcommand{\GG}{\mathbb{G}}

\newcommand{\NN}{\mathbb{Z}_{>0}}
\newcommand{\Zgezero}{\mathbb{Z}_{\geqslant 0}}

\newcommand{\kkk}{k}
\renewcommand{\lll}{l}
\newcommand{\uuu}{u}
\newcommand{\vvv}{v}
\newcommand{\lstar}{\,{}^c\!*}
\newcommand{\rstar}{*^c}
\newcommand{\wstar}{\,\widetilde*\,}
\newcommand{\lwstar}{\,{}^c\,\widetilde*\,}
\newcommand{\rwstar}{\,\widetilde*\,^{c}}

\DeclareMathOperator{\GL}{GL}
\DeclareMathOperator{\SL}{SL}

\DeclareMathOperator{\PSL}{PSL}

\DeclareMathOperator{\rk}{rk}

\theoremstyle{plain}

\newtheorem{proposition}{Proposition}
\newtheorem{theorem}{Theorem}

\theoremstyle{definition}

\newtheorem{example}{Example}

\newtheorem{question}{Question}
\theoremstyle{remark}
\newtheorem{remark}{Remark}

\newcommand\das{\renewcommand{\arraystretch}{1}}

\begin{document}

\title[Algebraic monoid structures on the affine 3-space]{Algebraic monoid structures on the affine 3-space}
\author{Ivan Arzhantsev}
\address{%
{\bf Ivan Arzhantsev} \newline\indent HSE University, Faculty of Computer Science, Pokrovsky Boulevard 11, Moscow, 109028 Russia}
\email{arjantsev@hse.ru}
\author{Roman Avdeev}
\address{%
{\bf Roman Avdeev} \newline\indent HSE University, Faculty of Computer Science, Pokrovsky Boulevard 11, Moscow, 109028 Russia}
\email{suselr@yandex.ru}
\author{Yulia Zaitseva}
\address{%
{\bf Yulia Zaitseva} \newline\indent HSE University, Faculty of Computer Science, Pokrovsky Boulevard 11, Moscow, 109028 Russia}
\email{yuliazaitseva@gmail.com}

\thanks{Supported by the RSF-DST grant 22-41-02019}



\begin{abstract}
We complete the classification of algebraic monoid structures on the affine 3-space. The result is based on a reduction of the general case to that of commutative monoids. 
We also study various algebraic properties of all monoids appearing in the classification.
\end{abstract}

\maketitle


\section{Introduction}

The theory of algebraic groups is one of the most rich and developed parts of mathematics that has deep connections with algebraic geometry, representation theory, number theory, topology, and combinatorics. The notion of algebraic monoid is a natural generalization of that of algebraic group; see~\cite{Pu1988,Re2005,Ri1,Vi} for general presentations. For algebraic monoids, we do not require that all elements be invertible, but the set of invertible elements of any algebraic monoid turns out to be an algebraic group and plays a significant role in studying the monoid. It is noted in the introduction to~\cite{Br-1} that ``from the viewpoint of algebraic geometry, it is an attractive problem to describe all algebraic semigroup structures on a given variety''.
In this paper, we focus on algebraic monoids whose underlying variety is the affine 3-space. 

Consider the affine space $\AA^n$ of dimension~$n$ over an algebraically closed field $\KK$ of characteristic zero.
An (\textit{algebraic}) \textit{monoid} structure on $\AA^n$ is a polynomial map $\mu \colon \AA^n \times \AA^n \to \AA^n$, $(x,y) \mapsto x * y$, satisfying the following properties:

\begin{enumerate}[label=\textup{(M\arabic*)},ref=\textup{M\arabic*}]
\item \label{M1}
$(x * y) * z = x * (y * z)$ for all $x,y,z \in \AA^n$;
\item \label{M2}
there is an element ${\bf 1} \in \AA^n$ such that $x * {\bf 1} = x = {\bf 1} * x$ for all $x \in \AA^n$.
\end{enumerate}

\noindent In addition, a monoid structure $\mu$ on $\AA^n$ is said to be \textit{commutative} if $x * y = y * x$ for all $x, y \in \AA^n$.

Two algebraic monoid structures $\mu_1$ and $\mu_2$ on $\AA^n$ are called \textit{isomorphic} if there is a polynomial automorphism $\varphi \colon \AA^n \to \AA^n$ such that $\varphi (x *_1 y) = \varphi(x) *_2 \varphi(y)$ for all $x,y \in \AA^n$.

A natural problem is to classify all monoid structures on $\AA^n$ up to isomorphism.
The case $n = 1$ is quite elementary and was considered in~\cite{Yo63}.
It turns out that there are only two non-isomorphic monoid structures on $\AA^1$, given by
\begin{equation} \label{eqn_monA1}
x * y = x + y \quad \text{and} \quad x * y = x \cdot y.
\end{equation}
Indeed, such a structure is defined by a single polynomial $f(x,y)$ in two variables, for which property (\ref{M1}) translates into the polynomial identity
\begin{equation} \label{A1assoc}
f(f(x,y),z) = f(x,f(y,z)).
\end{equation}
It implies that each monomial in $f$ contains each of the variables $x,y$ with the power at most~$1$.
Thus $f$ is a linear combination of the monomials $xy,x,y,1$, and condition~(\ref{A1assoc}) is equivalent to several relations between the coefficients.
Adding property (\ref{M2}) and taking into account that all polynomial automorphisms of $\AA^1$ have the form $x \mapsto ax+b$ with $a,b\in \KK$ and $a \ne 0$, one obtains~(\ref{eqn_monA1}).

\smallskip

The case $n=2$ is more complicated; all commutative (resp. noncommutative) monoid structures on~$\AA^2$ were classified in~\cite{ABZ2020} (resp.~\cite{Bi2021}).
As a result, all pairwise non-isomorphic monoid structures on $\AA^2$ are given by the following formulas: 

\smallskip

$(x_1,y_1)*(x_2,y_2) = (x_1+x_2, y_1 + y_2)$;

\smallskip

$(x_1,y_1)*(x_2,y_2) = (x_1x_2,x_1^ay_2+x_2^by_1)$ ($a,b \in \ZZ_{\ge0}$);

\smallskip

$(x_1,y_1)*(x_2,y_2) = (x_1x_2, y_1y_2)$.

\smallskip

For the case $n=3$, several partial classification results are known from~\cite{ABZ2020} and~\cite{Za2024}; see details in Section~\ref{prelim_sec}. These include the classification of all commutative monoid structures on~$\AA^3$ obtained in~\cite{ABZ2020}.

In the present paper we complete the classification of all monoid structures on~$\AA^3$.
The main idea consists in reducing the noncommutative case to the commutative one.

To state our result, we need some additional notation. 
A quadruple $p  = (b,b',c,c') \in \NN^4$ will be called \emph{compatible} if there exists an integer $d > 0$ such that
\begin{equation}
\label{bcbcd_cond_eq}
\begin{aligned}
c &- bd \ge 0, \quad\;\, c'-b'd \ge 0,\\
c &-b(d+1) = c'-b'(d+1) < 0.
\end{aligned}
\end{equation}
Here $d$ is the quotient of $c$ modulo~$b$ and also of $c'$ modulo~$b'$. 
If $b \ne b'$, then $p$ is compatible if and only if the numbers $b,b',c,c'$ satisfy the condition
\begin{equation}
\label{bcbc_cond_eq}
\biggl\lfloor\frac{c}{b}\biggr\rfloor = \biggl\lfloor\frac{c'}{b'}\biggr\rfloor = \frac{c-c'}{b-b'} - 1 > 0.
\end{equation}
If $b = b'$, then $p$ is compatible if and only if $c = c'$ and $c \ge b$.
For a compatible quadruple~$p$, we define the polynomial
\begin{multline*}
Q_p(x_1, x_2, y_1, y_2) 
= \sum_{\substack{i + j = d+1 \\ i,j \ge 1}} \binom{d+1}{i}x_1^{c-bi}x_2^{c'-b'j}y_1^i y_2^j = 
\\
= x_1^cx_2^{c'}\left(\biggl(\frac{y_1}{x_1^b}+\frac{y_2}{x_2^{b'}}\biggr)^{\!\!d+1} \!\!- \biggl(\frac{y_1}{x_1^b}\biggr)^{\!\!d+1} \!\!- \biggl(\frac{y_2}{x_2^{b'}}\biggr)^{\!\!d+1}\right).
\end{multline*}

Let us call quadruples $(b,b',c,c')$ and $(c,c',b,b')$ \emph{symmetric} to each other.

\smallskip

The main result of this paper is given by the following theorem.

\renewcommand\arraystretch{1.5}

\begin{theorem} \label{A3_theor}
The following assertions hold.
\begin{enumerate}[label=\textup{(\alph*)},ref=\textup{\alph*}]
\item \label{A3_theor_a}
Up to isomorphism, every monoid structure on $\AA^3$ appears in Table~\textup{\ref{table_ams_on_A3}}, where we use the notation $p = (b, b', c, c')$. 

\item \label{A3_theor_b}
In Table~\textup{\ref{table_ams_on_A3}}, any two monoid structures of different types are not isomorphic; two monoid structures of the same type with different values of parameters are isomorphic if and only if they are of type $\mathrm{1M2A}_0(p)$ or $\mathrm{2M1A}(p)$ and the corresponding quadruples $p = (b,b',c,c')$ are symmetric to each other.
\end{enumerate}
\end{theorem}

\begin{table}[h]
 \centering
\caption{Algebraic monoid structures on~$\AA^3$}
\label{table_ams_on_A3}
\(
\begin{array}{|c|c|c|}
\hline \text{rk} & \text{Notation} & 
(x_1, y_1, z_1) * (x_2, y_2, z_2) \\\hline\hline
0 & \mathrm{3A} &
(x_1 + x_2, \, y_1 + y_2, \, z_1 + z_2)
\\\hline
0 & \mathrm{U_3} &
(x_1 + x_2, \, y_1 + y_2, \, z_1 + z_2 + x_1y_2)
\\\hline
1 & \text{1M2A}_0(p) &
\bigl(x_1x_2, \, x_1^b y_2 + x_2^{b'} y_1, \, x_1^c z_2 + x_2^{c'} z_1\bigr)\;
\text{ for } p \in \Zgezero^4
\\\hline
1 & \text{1M2A}_1(p) &
  \begin{matrix}
  \bigl(x_1x_2, \, x_1^b y_2 + x_2^{b'} y_1, \, x_1^c z_2 + x_2^{c'} z_1 +
  Q_p(x_1, x_2, y_1, y_2)\bigr) \\
  \text{ for compatible } p \in \NN^4
  \end{matrix}
\\\hline
2 & \text{2M1A}(p) &
(x_1 x_2, \, y_1 y_2, \, x_1^by_1^c z_2 + x_2^{b'} y_2^{c'} z_1) \; 
\text{ for } p \in \Zgezero^4
\\\hline
3 & \mathrm{3M} &
(x_1 x_2, \, y_1 y_2, \, z_1 z_2)\\\hline
\end{array}
\)
\end{table}

In the first column of Table~\ref{table_ams_on_A3}, for each monoid structure we indicate the value of its rank; in the second column, the notation $k\mathrm{M}l\mathrm{A}$ reflects that the group of invertible elements of the corresponding monoid is isomorphic to a semidirect product $\GG_m^k \rightthreetimes \GG_a^l$; see Section~\ref{prelim_sec} for the definitions of rank and group of invertible elements.

One can see that the only groups in Theorem~\ref{A3_theor} are $\mathrm{3A}$ and~$\mathrm{U_3}$, which are the $3$-dimensional vector group and the group of upper unitriangular $3\times3$-matrices, respectively.

It is worth mentioning that, as a consequence of Lazard's theorem~\cite{Laz55} (see also~\cite[Remark~3 and Theorem~11]{Pop22}), all algebraic groups whose underlying variety is an affine space (of arbitrary dimension) are precisely the unipotent ones.

Let us illustrate Theorem~\ref{A3_theor} by example. The quadruple $p = (b,b',c,c') = (1,2,1,3)$ is compatible since $\left\lfloor\frac{1}{1}\right\rfloor = \left\lfloor\frac{3}{2}\right\rfloor = \frac{1-3}{1-2} - 1 > 0$. Then the formula
\[(x_1, y_1, z_1) * (x_2, y_2, z_2) = \bigl(x_1x_2, \, x_1 y_2 + x_2^2 y_1, \, x_1 z_2 + x_2^3 z_1 + 2 x_2y_1y_2\bigr)\]
defines a noncommutative monoid structure on $\AA^3$ of type $\mathrm{1M2A}_1(p)$. 

This paper is organized as follows. In Section~\ref{prelim_sec} we collect and discuss all notions and results on algebraic monoids needed in the paper.
In Section~\ref{reduction_sec} we present the construction of a reduction to the commutative case, which is a key ingredient in the proof of Theorem~\ref{A3_theor}.
The proof itself is carried out in Section~\ref{proof_sec}. In Section~\ref{properties_sec} we study various algebraic properties of all monoid structures on~$\AA^3$ including the structure of the center, set of idempotents, kernel, and existence of a zero. We also pose several questions on these objects for monoid structures on affine spaces of higher dimensions.
Finally, in the Appendix, which is independent of the main part of the paper, we classify all reductive monoid structures (see their definition in Section~\ref{prelim_sec}) on affine spaces of arbitrary dimension.

\smallskip

\textbf{Acknowledgement.}
We thank the referee for useful comments on a previous version of this paper.

\section{Preliminaries}
\label{prelim_sec}

We work over an algebraically closed field~$\KK$ of characteristic zero. The symbol $\GG_m$ (resp.~$\GG_a$) stands for the multiplicative (resp. additive) group of~$\KK$ regarded as an algebraic group.

Let $X$ be an algebraic variety (possibly reducible). An (\emph{algebraic}) \emph{monoid structure} on~$X$ is a morphism $\mu\colon X \times X \to X$, $\mu(x,y) = x*y$, satisfying the \emph{associativity} property $(x*y)*z = x*(y*z)$ for all $x,y,z \in X$ and possessing an element~${\bf 1} \in X$, called a \emph{unity}, such that ${\bf 1} * x = x * {\bf 1} = x$ for all $x \in X$.
In this situation, the pair $(X,\mu)$ is said to be an (\emph{algebraic}) \emph{monoid} with underlying variety~$X$ and multiplication~$\mu$.

A monoid structure $\mu$ on~$X$ (and also the monoid $(X,\mu)$ itself) is called \emph{commutative} if $x * y = y * x$ for all $x,y \in X$.

Two monoids $(X,\mu)$ and $(\widetilde X, \widetilde \mu)$ are said to be \textit{isomorphic} if there is an isomorphism of algebraic varieties $\varphi \colon X \xrightarrow{\sim} \widetilde X$ such that $\widetilde\mu(\varphi(x),\varphi(y)) = \varphi(\mu(x,y))$ for all $x,y \in X$.
If moreover $X = \widetilde X$, then we simply say that $\mu$ and $\widetilde \mu$ are isomorphic monoid structures on~$X$.

For every monoid $(X,\mu)$, the set of invertible elements in~$X$ is Zariski open and thus an algebraic group; see~\cite[Theorem~1]{Ri1} and~\cite[Theorem~1]{Br-1}. We denote this group by $G(X,\mu)$.
Moreover, it is known that $X$ is affine if and only if $G(X,\mu)$ is a linear algebraic group; see~\cite[Theorem~5]{Ri2}. 

Given monoids $(X_1,\mu_1), (X_2,\mu_2), \ldots, (X_k,\mu_k)$, one defines their direct product $(X,\mu) = (X_1,\mu_1) \times (X_2,\mu_2) \times \ldots \times (X_k,\mu_k)$ where $X = X_1\times X_2 \times \ldots \times X_k$ and the multiplication~$\mu$ is given by
\[
(x_1, x_2, \ldots, x_k) * (y_1, y_2, \ldots y_k) = (\mu_1(x_1,y_1), \mu_2(x_2,y_2), \ldots, \mu_k(x_k,y_k)).
\]
Clearly, $(X,\mu)$ is a monoid and $G(X,\mu) = G(X_1,\mu_1)\times \ldots \times G(X_k,\mu_k)$.

Given an algebraic group~$G$, a \emph{group embedding} of~$G$ is an open embedding $G \hookrightarrow Y$ into an algebraic variety~$Y$ such that the action of the group $G \times G$ on $G$ defined by $(g_1,g_2)\cdot g = g_1gg_2^{-1}$ can be extended to an action of $G \times G$ on~$Y$.

Let $X$ be an irreducible affine variety. For every monoid structure $\mu$ on~$X$, the natural inclusion $G(X,\mu) \hookrightarrow X$ is clearly a group embedding, where the required action of $G(X,\mu) \times G(X,\mu)$ on~$X$ comes from the monoid multiplication: $(g_1, g_2) \cdot x = g_1 * x * g_2^{-1}$. In fact, the converse is also true: for every group embedding $G \hookrightarrow X$ of a connected linear algebraic group~$G$, the multiplication map $G \times G \to G$ extends to a monoid structure on~$X$; this follows from~\cite[formula~(3)]{Vi} in the case of reductive~$G$ in characteristic zero and was proved in~\cite[Proposition~1]{Ri1} in the general case (including arbitrary characteristic).
In this paper, we will need the following consequence of the above correspondence.

\begin{proposition} \label{comm_mon}
Let $K$ be a connected commutative algebraic group and let $K \hookrightarrow Y$ be an open embedding into an irreducible affine variety~$Y$. Suppose that the action of $K$ on itself by left multiplication, $(g,y) \mapsto gy$, extends to an action of $K$ on the whole~$Y$. Then the action morphism $K \times Y \to Y$ extends to a morphism $Y \times Y \to Y$, which is a commutative monoid structure on $Y$ with set of invertible elements~$K$.
\end{proposition}

\begin{proof}
As $K$ is commutative, the action of $K$ on itself by right multiplication, $(g,y) \mapsto yg^{-1}$, also extends to an action on the whole~$Y$. Thus $K \hookrightarrow Y$ is a group embedding and the above correspondence yields a monoid structure $\mu \colon Y \times Y \to Y$ extending the action morphism $K \times Y \to Y$.
Since $\mu$ is commutative on the open subset $K \times K$, it is commutative on the whole $X \times X$.
\end{proof}

In the remaining part of this section we present known classification results for affine algebraic monoids.
We assume that $X$ is a normal irreducible affine variety and put $n =\dim X$. Given a monoid structure $\mu$ on~$X$, the \emph{rank} of $(X,\mu)$ is the dimension of a maximal torus in $G(X,\mu)$; we denote it by~$\rk(X,\mu)$.

A monoid $(X,\mu)$ is called \emph{reductive} if the group $G(X,\mu)$ is reductive. A combinatorial description of all reductive monoids is obtained in~\cite{Vi,Ri1}.

Among well-understood monoids $(X,\mu)$ are those with $\rk(X,\mu) \in \lbrace 0,n-1,n\rbrace$.
If $\rk (X,\mu) = 0$, then the group $G(X,\mu)$ is unipotent.
By~\cite[Section~1.3]{PV}, all orbits of unipotent group actions on affine varieties are closed, therefore $G(X,\mu) = X$. So in this case all monoids $(X,\mu)$ are precisely all unipotent groups; in particular, $X \simeq \AA^n$.
If $\rk(X,\mu) = n$, then the group $G(X,\mu) = T$ is a torus and thus $X$ is a group embedding of~$T$ (in particular, $(X,\mu)$ is automatically commutative). Such varieties are called \emph{toric varieties}; moreover, the toric structure on them is unique up to isomorphism; see~\cite{BB} for $X = \AA^n$, \cite{De1982,Gu1998} for the affine case, and~\cite{Be2003} for the general case. So the classification of all monoids $(X,\mu)$ with $\rk(X,\mu) = n$ is equivalent to the classification of all affine toric varieties, which is well known and given in terms of cones; see~\cite[Section~1.3]{Fu1993} or~\cite[Chapter~1]{CLS2011}.
All monoids $(X,\mu)$ with $\rk(X,\mu) = n-1$ were classified in~\cite{ABZ2020} for $X = \AA^n$ and commutative~$\mu$, in~\cite{DzZa2021} for arbitrary $X$ and commutative~$\mu$, in~\cite{Bi2021} for $n=2$ and noncommutative~$\mu$, and in~\cite{Za2024} for arbitrary~$X$ and noncommutative~$\mu$.
We remark that the results mentioned in this paragraph imply the classification of all monoids $(X,\mu)$ with $n\le2$.

All commutative monoids $(X,\mu)$ with $X = \AA^3$ and $\rk(X,\mu)=1$ were classified in~\cite{ABZ2020}. Combining this with the above-mentioned results yields the classification of all commutative monoid structures on~$\AA^3$, which is given in~\cite[Theorem~1]{ABZ2020}. This classification is essentially used in this paper and reproduced in Theorem~\ref{A3comm_theor} below. It involves polynomials $Q_{b,c}$ with $b, c \in \NN$, which can be defined as follows:
\begin{multline*}
Q_{b,c}(x_1, x_2, y_1, y_2) =
\sum \limits_{\substack{i + j = d+1 \\ i,j \ge 1}} \binom{d+1}{i} x_1^{c-bi} x_2^{c-bj} y_1^i y_2^j = \\ 
= \frac{(x_1^by_2+x_2^by_1)^{d+1}-(x_1^by_2)^{d+1}-(x_2^by_1)^{d+1}}
{x_1^{b-e} y_1^{b-e}},
\end{multline*}
where $c = b d + e; \; d, e \in \ZZ; \; 0 \le e < b$. 

\renewcommand\arraystretch{1.5}

\begin{theorem}[{\cite[Theorem~1]{ABZ2020}}] \label{A3comm_theor}
The following assertions hold.
\begin{enumerate}[label=\textup{(\alph*)},ref=\textup{\alph*}]
\item \label{A3comm_theor_a}
Up to isomorphism, all commutative monoid structures on $\AA^3$ appear in Table~\textup{\ref{table_cams_on_A3}}, where we use the notation from Table~\textup{\ref{table_ams_on_A3}}.

\item \label{A3comm_theor_b}
In Table~\textup{\ref{table_cams_on_A3}}, any two monoid structures of different types are not isomorphic; the only isomorphisms between two monoid structures of the same type with different values of parameters are $\mathrm{1M2A}_0(b,b,c,c) \simeq \mathrm{1M2A}_0(c,c,b,b)$ and $\mathrm{2M1A}(b,b,c,c) \simeq \mathrm{2M1A}(c,c,b,b)$.
\end{enumerate}
\end{theorem}

\begin{table}[h]
 \centering
\caption{Commutative monoid structures on~$\AA^3$}
\label{table_cams_on_A3}
\(
\begin{array}{|c|c|c|}\hline
\text{rk} & \text{Monoid} & (x_1, y_1, z_1) * (x_2, y_2, z_2) \\\hline\hline
0 & \mathrm{3A} &
(x_1 + x_2, \, y_1 + y_2, \, z_1 + z_2)
\\\hline
1 & \mathrm{1M2A}_0(b,b,c,c) &
\bigl(x_1x_2, \, x_1^b y_2 + x_2^b y_1, \, x_1^c z_2 + x_2^c z_1\bigr),
\; b, c \in \Zgezero
\\\hline
1 & \mathrm{1M2A}_1(b,b,c,c) &
\begin{matrix}
\bigl(x_1x_2, \, x_1^b y_2 + x_2^b y_1, \, x_1^c z_2 + x_2^c z_1 + Q_{b, c}(x_1, x_2, y_1, y_2)\bigr), \\
b, c \in \NN, \; b \le c
\end{matrix}
\\\hline
2 & \mathrm{2M1A}(b,b,c,c) &
(x_1 x_2, \, y_1 y_2, \, x_1^by_1^c z_2 + x_2^b y_2^c z_1), \; b, c \in \Zgezero
\\\hline
3 & \mathrm{3M} &
(x_1 x_2, \, y_1 y_2, \, z_1 z_2)\\\hline
\end{array}
\)
\end{table}

We remark that in~\cite{GvSe2022} Theorem~\ref{A3comm_theor} was extended to the case of non-closed ground fields of characteristic zero.

Observe that isomorphisms $\mathrm{1M2A}_0(b,b,c,c) \simeq \mathrm{1M2A}_0(c,c,b,b)$ and $\mathrm{2M1A}(b,b,c,c) \simeq \mathrm{2M1A}(c,c,b,b)$ are given by $(x,y,z) \mapsto (x,z,y)$ and $(x,y,z) \mapsto (y,x,z)$, respectively.

\medskip

As follows from the discussion in this section, in the classification of all monoid structures on~$\AA^3$ the only remaining case not covered by known results is that of noncommutative monoid structures of rank~$1$, and we complete it in this paper.

\section{Reduction to the commutative case}
\label{reduction_sec}

In this section, for a not necessarily commutative monoid $(X,\mu)$, we construct two commutative monoid structures on $X$ in the case where the group $G(X,\mu)$ is isomorphic to a semidirect product of a torus $T$ and a commutative unipotent group $U$.
The main idea is to consider the action on~$X$ of the direct product $K = T \times U$ where $T$ (resp.~$U$) acts by left (resp. right) multiplication. This action satisfies the conditions of Proposition~\ref{comm_mon} and thus yields a commutative monoid structure on $X$ with group of invertible elements~$K$. The second monoid structure is obtained similarly by involving the right multiplication by $T$ and left multiplication by~$U$. Below we provide a more general construction, which works for a semidirect product of any two commutative linear algebraic groups.

Let $X$ be an irreducible affine variety.
We say that a monoid structure $\mu$ on~$X$ (and also the corresponding monoid $(X,\mu)$) is \textit{semicommutative} if the group $G(X,\mu)$ is isomorphic to a semidirect product of two commutative linear algebraic groups (which are necessarily connected since $X$ is irreducible).

In what follows, we assume that $\mu$ is a semicommutative monoid structure on~$X$. For convenience, we use the notation $G = G(X,\mu)$.
Let $X_0 \subseteq X$ be the set of invertible elements of~$X$.
We distinguish between the notations $G$ and $X_0$ since later we will consider algebraic group structures on~$X_0$ different from~$G$.
Fix commutative subgroups $T,U \subseteq G$ such that $G = T \rightthreetimes U$.
Fix a group isomorphism $X_0 \simeq G$ and regard $T,U$ as subvarieties of~$X_0$ via this isomorphism.
Let ${\bf 1} \in X_0$ be the unity of $(X, \mu)$. Consider the commutative group $K = T \times U$.

Since $G = T U = U T$, any element $g \in X_0$ can be uniquely represented in two ways: 
\[g =t * \uuu = \vvv * t, \quad\text{ where } \uuu,\vvv \in U \subseteq X_0, \; t \in T \subseteq X_0.\]
Consider the following two morphisms:
\begin{equation} \label{r_op_eq}
X_0 \times X \to X, \quad (t*\uuu, \, x) \mapsto t*x*\uuu;\\
\end{equation}
\begin{equation} \label{l_op_eq}
X_0 \times X \to X, \quad (\vvv*t, \, x) \mapsto \vvv*x*t.
\end{equation}

\begin{proposition} \label{comm_red_prop}
Morphism~\eqref{r_op_eq} extends to a morphism
\begin{equation} \label{cmu}
{}^c\!\mu\colon X \times X \to X, \quad (x, y) \mapsto x\,^c\!*y, 
\end{equation}
and morphism~\eqref{l_op_eq} extends to a morphism
\begin{equation} \label{muc}
\mu^c\colon X \times X \to X, \quad (x, y) \mapsto x*^cy. 
\end{equation}
Both morphisms define commutative monoid structures on $X$ with set of invertible elements~$X_0$ isomorphic to~$K$ as an algebraic group.
\end{proposition}

\begin{proof}
Consider the action of $K$ on~$X$ given by $((t,u),x) \mapsto t*x*u$. Then the orbit of~${\bf 1}$ is~$X_0$ and the stabilizer of~${\bf 1}$ is trivial, so the map $K \to X$ given by $(t,u) \mapsto t*u$ is an open embedding. Identifying~$K$ with its image~$X_0$ we see that morphism~\eqref{r_op_eq} extends the action of~$K$ on itself by left multiplication. By Proposition~\ref{comm_mon}, morphism~\eqref{r_op_eq} extends further to a morphism ${}^c\!\mu \colon X \times X \to X$, which yields a commutative monoid structure on~$X$ with set of invertible elements~$X_0$ isomorphic to~$K$ as an algebraic group.

The result for~\eqref{l_op_eq} is obtained by applying the same argument to the action of $K$ on~$X$ given by  $((t,v),x) \mapsto v*x*t$.
\end{proof}

The monoid structure ${}^c\!\mu$ given by~\eqref{cmu} (resp. the monoid $(X,{}^c\!\mu)$) is called the \emph{left commutative reduction} of $\mu$ (resp.~$(X,\mu)$). Similarly, the monoid structure $\mu^c$ given by~\eqref{muc} (resp. the monoid $(X,\mu^c)$) is called the \emph{right commutative reduction} of $\mu$ (resp.~$(X,\mu)$). It follows from the construction that in general ${}^c\!\mu$ and $\mu^c$ depend on the choices of~$T$ and~$U$ as subvarieties of~$X_0$.

In the following two examples the fact that $\mu$ is a monoid structure on~$X$ follows from the proof of Theorem~\ref{A3_theor}(\ref{A3_theor_a}) given in Section~\ref{proof_sec}. These examples will be used in the proof of Theorem~\ref{A3_theor}(\ref{A3_theor_b}).

\begin{example}
\label{MMA_red}
Take $X=\AA^3$ with coordinates $x,y,z$ and consider a monoid $(X,\mu)$ of type $\mathrm{2M1A}(p)$ for some $p = (b,b',c,c') \in \Zgezero^4$, so that the multiplication is given by
\[(x_1, y_1, z_1) * (x_2, y_2, z_2) = (x_1 x_2, \, y_1 y_2, \, x_1^by_1^c z_2 + x_2^{b'} y_2^{c'} z_1).\]
The set of invertible elements is $X_0= \{x,y \ne 0\}$; as an algebraic group it is a semidirect product of $T = \{x,y \ne 0; \ z=0\} \simeq \GG_m^2$ and $U = \{x=y=1\} \simeq \GG_a$, so $(X,\mu)$ is semicommutative. Let us calculate the left commutative reduction of $(X,\mu)$ with respect to $T$ and~$U$.
An element $g = (x,y,z) \in X_0 \subseteq \AA^3$ is decomposed as $g = t*\uuu = \vvv*t$, where
\[t = (x,y,0)\in T, \;\; \uuu = (1,1,x^{-b}y^{-c}z) \in U, \;\; \vvv = (1,1,x^{-b'}y^{-c'}z) \in U.\]
Morphism~\eqref{r_op_eq} is given by
\begin{multline*}
(x,y,z) \,{}^{c}\!* (x_2, y_2, z_2) = t * (x_2, y_2, z_2) * \uuu = 
\\
= (xx_2, yy_2, x^by^cz_2) * (1,1,x^{-b}y^{-c}z) = (xx_2, yy_2, x^by^cz_2 + x_2^by_2^cz),
\end{multline*}
so the left commutative reduction of $\mathrm{2M1A}(p)$ is $\mathrm{2M1A}(b,b,c,c)$. Similar calculations show that the right commutative reduction is $\mathrm{2M1A}(b',b',c',c')$.
\end{example}

\begin{example}
\label{MAA_red}
Take again $X=\AA^3$ with coordinates $x,y,z$ and consider a monoid $(X,\mu)$ of type $\mathrm{1M2A}_\kappa(p)$ for some $\kappa \in \lbrace 0,1 \rbrace$ and $p = (b,b',c,c') \in \Zgezero^4$, so that the multiplication is given by
\[(x_1, y_1, z_1) * (x_2, y_2, z_2) = \bigl(x_1x_2, \, x_1^b y_2 + x_2^{b'} y_1, \, x_1^c z_2 + x_2^{c'} z_1 + \kappa Q_p(x_1, x_2, y_1, y_2)\bigr).\]
Then $X_0 = \{x \ne 0\}$; as an algebraic group it is a semidirect product of $T = \{x \ne 0; \; y,z=0\} \simeq \GG_m$ and $U = \{x=1\} \simeq \GG_a^2$, so $(X,\mu)$ is semicommutative. Let us calculate the right commutative reduction of $(X,\mu)$ with respect to $T$ and~$U$.
An element $g=(x,y,z) \in X_0$ is decomposed as $g=t*\uuu = \vvv*t$, where
\[t = (x,0,0)\in T, \;\; \uuu = (1,x^{-b}y,x^{-c}z) \in U, \;\; \vvv = (1,x^{-b'}y, x^{-c'}z) \in U.\]
Morphism~\eqref{l_op_eq} is given by
\begin{multline*}
(x_1,y_1,z_1) *^c (x, y, z) = \vvv * (x_1, y_1, z_1) * t = 
(1,x^{-b'}y, x^{-c'}z) * (x_1x, x^{b'}y_1, x^{c'}z_1) = 
\\
= \bigl(x_1x, \, x_1^{b'} y + x^{b'} y_1, \, x_1^{c'} z + x^{c'} z_1 + \kappa Q_p(1, x_1x, x^{-b'}y, x^{b'}y_1)\bigr).
\end{multline*}
Since $Q_p(1, x_1x, x^{-b'}y, x^{b'}y_1) = Q_{b',c'}(x_1, x, y_1, y)$, we conclude that the right commutative reduction of $\mathrm{1M2A}_0(p)$ is $\mathrm{1M2A}_0(b',b',c',c')$ and that of $\mathrm{1M2A}_1(p)$ is $\mathrm{1M2A}_1(b',b',c',c')$. In the same way, the left commutative reductions are $\mathrm{1M2A}_0(b,b,c,c)$ and $\mathrm{1M2A}_1(b,b,c,c)$, respectively. 
\end{example}

Let $(\widetilde X, \widetilde \mu)$ be one more semicommutative monoid with set of invertible elements $\widetilde X_0$ and commutative subgroups $\widetilde T, \widetilde U \subseteq \widetilde X_0$ such that $\widetilde X_0 = \widetilde T \rightthreetimes \widetilde U$ as algebraic groups.
Let $(\widetilde X, {}^c\widetilde\mu)$ and $(\widetilde X, \widetilde\mu^c)$ be the left and right commutative reductions of $(\widetilde X, \widetilde \mu)$, respectively, with respect to~$\widetilde T$ and~$\widetilde U$.

\begin{proposition}
\label{muccmu_prop}
Suppose that $\varphi \colon X \to \widetilde X$ is an isomorphism of algebraic varieties such that $\varphi(T) = \widetilde T$ and $\varphi(U) = \widetilde U$. Then the following conditions are equivalent:
\begin{enumerate}[label=\textup{(\arabic*)},ref=\textup{\arabic*}]
\item \label{mon_isom_a}
$\varphi$ induces a monoid isomorphism $(X,\mu) \xrightarrow{\sim} (\widetilde X, \widetilde\mu)$;

\item \label{mon_isom_b}
$\varphi$ simultaneously induces monoid isomorphisms $(X, {}^c\!\mu) \xrightarrow{\sim} (\widetilde X, {}^c\widetilde \mu)$ and $(X,\mu^c) \xrightarrow{\sim} (\widetilde X, \widetilde\mu^c)$.
\end{enumerate}
\end{proposition}

\begin{proof}
First of all, note that $\varphi$ induces a monoid isomorphism $(X,\mu) \xrightarrow{\sim} (\widetilde X, \widetilde\mu)$, $(X, {}^c\!\mu) \xrightarrow{\sim} (\widetilde X, {}^c\widetilde \mu)$, or $(X,\mu^c) \xrightarrow{\sim} (\widetilde X, \widetilde\mu^c)$
if and only if the equality
\begin{gather}
\varphi(x * y) = \varphi(x) \wstar \varphi(y), \label{mon_isom1}\\
\varphi(x \lstar y) = \varphi(x) \lwstar \varphi(y), \label{mon_isom2} \\
\varphi(x \rstar y) = \varphi(x) \rwstar \varphi(y), \label{mon_isom3}
\end{gather}
respectively, holds for all $x,y \in X$.

(\ref{mon_isom_a})$\Rightarrow$(\ref{mon_isom_b})
Take any $t \in T$, $u \in U$, $y \in X$ and put $x = t*u$. Using~(\ref{mon_isom1}) we obtain
\[
\varphi(x \lstar y) = \varphi (t * y * u) = \varphi(t) \wstar \varphi(y) \wstar \varphi(u) = (\varphi(t) \wstar \varphi(u)) \lwstar \varphi(y) = \varphi(x) \lwstar \varphi(y),
\]
which yields~(\ref{mon_isom2}) for all $x \in X_0$ and $y \in X$. Since $X_0 \times X$ is open in~$X \times X$, it follows that~(\ref{mon_isom2}) holds for all $x,y \in X$.
In a similar way, one proves~(\ref{mon_isom3}).

(\ref{mon_isom_b})$\Rightarrow$(\ref{mon_isom_a})
As above, it suffices to show~(\ref{mon_isom1}) for all $x \in X_0$ and $y \in X$.
Taking $t \in T$, $u \in U$ such that $x = t*u$ and using~(\ref{mon_isom2}),\,(\ref{mon_isom3}) we obtain
\begin{multline*}
\varphi(x * y) = \varphi((t * u) * y) = \varphi(t * (u * y)) = \varphi (t \lstar (u * y)) = \varphi (t) \lwstar \varphi(u * y) = \\
\varphi (t) \wstar \varphi(u * y) = \varphi (t) \wstar \varphi(u \rstar y) = \varphi (t) \wstar (\varphi(u) \rwstar \varphi(y)) = \varphi (t) \wstar (\varphi(u) \wstar \varphi(y)) = \\
(\varphi (t) \wstar \varphi(u)) \wstar \varphi(y) = (\varphi (t) \lwstar \varphi(u)) \wstar \varphi(y) = \varphi (t \lstar u) \wstar \varphi(y) = \\
\varphi (t * u) \wstar \varphi(y) = \varphi(x) \wstar \varphi(y)
\end{multline*}
as required.
\end{proof}

\begin{proposition}
\label{choiceTU_prop}
Suppose that $T$ is a torus and $U$ is a commutative unipotent group.
Then the isomorphism classes of left and right commutative reductions of~$X$ do not depend on the choices of $T$ and~$U$ as subvarieties of~$X_0$.
\end{proposition}

\begin{proof}
Clearly, $U$ is the unipotent radical of~$G$, so it is identified with a subvariety of~$X_0$ in a unique way. Next, $T$ is a maximal torus of~$G$, hence it is uniquely determined up to conjugation. Let $\widetilde T \subseteq X_0$ be another choice of a maximal torus. If $a \in X_0$ is an element such that $a*T*a^{-1} = \widetilde T$, then the morphism $X \to X$, $x \mapsto a*x*a^{-1}$, induces an automorphism of the monoid~$(X,\mu)$ that preserves~$U$ and maps $T$ to $\widetilde T$. Now the claim follows from Proposition~\ref{muccmu_prop}.
\end{proof}

\begin{remark}
In contrast to Proposition~\ref{muccmu_prop}, it may happen that two non-isomorphic monoid structures $\mu$ and $\widetilde \mu$ on the same variety~$X$ have isomorphic left commutative reductions and isomorphic right commutative reductions.
For example, consider the monoid structures $\mu$ of type $\mathrm{1M2A}_0(b,b,c,c)$ and $\widetilde \mu$ of type $\mathrm{1M2A}_0(b,c,c,b)$ on $X = \AA^3$ with $b \ne c$.
The first one is commutative and the second one is not, so they are not isomorphic. However, according to Example~\ref{MMA_red}, both ${}^c\!\mu$ and ${}^c\widetilde\mu$ are isomorphic to $\mathrm{1M2A}_0(b,b,c,c)$ and both $\mu^c$ and $\widetilde\mu^c$ are isomorphic to $\mathrm{1M2A}_0(b,b,c,c) \simeq \mathrm{1M2A}_0(c,c,b,b)$. In view of Proposition~\ref{muccmu_prop} this means that there exists no isomorphism of varieties $\varphi \colon X \xrightarrow{\sim} X$ inducing simultaneously monoid isomorphisms $(X, {}^c\!\mu) \xrightarrow{\sim} (X, {}^c\widetilde \mu)$ and $(X,\mu^c) \xrightarrow{\sim} (X, \widetilde\mu^c)$
\end{remark}

\begin{remark}
The construction of left (or right) commutative reduction provides an approach to classification of noncommutative monoids $(X,\mu)$ with $\rk(X,\mu) = \dim X - 1$ different to that used in~\cite{Za2024}.
\end{remark}

\section{Proof of Theorem~\ref{A3_theor}}
\label{proof_sec}

In this section we prove Theorem~\ref{A3_theor}.

The main idea of the proof of part~(\ref{A3_theor_a}) is as follows. For any monoid $(\AA^3,\mu)$ with set of invertible elements~$X_0$ and $G(X,\mu) \simeq \GG_m^p \rightthreetimes \GG_a^q$, $p+q=3$, we consider its left commutative reduction $(X,{}^c\!\mu)$. According to Proposition~\ref{comm_red_prop}, the set of invertible elements of $(X,{}^c\!\mu)$ is also $X_0$ and $G(X,{}^c\!\mu) \simeq K$ with $K = \GG_m^p \times \GG_a^q$, which yields a group embedding $K \hookrightarrow \AA^3$. All such embeddings were classified in~\cite{ABZ2020}, and we check for which of them the multiplication $\mu$ on the image of~$K$ in~$\AA^3$ extends to the whole~$\AA^3$.

In the proof of part~(\ref{A3_theor_b}) we consider both left and right commutative reductions and use Propositions~\ref{muccmu_prop} and~\ref{choiceTU_prop}.

In our arguments we will need that every semidirect product of the form $\GG_m^p \rightthreetimes \GG_a^q$ is defined by a set of characters $\chi_1,\ldots, \chi_q$ of $\GG_m^p$ in such a way that the adjoint action of $\GG_m^p$ on $\GG_a^q$ is by the formula $t \cdot (\alpha_1,\ldots, \alpha_q) = (\chi_1(t)\alpha_1,\ldots, \chi_q(t)\alpha_q)$. 
In particular, a semidirect product $\GG_m \rightthreetimes \GG_a^2$ is defined by two characters $t \mapsto t^{-\kkk}$ and $t \mapsto t^{-\lll}$ of $\GG_m$, where $\kkk,\lll \in \ZZ$. (The minus sign here is chosen for convenience in subsequent computations.) We denote this group by $G_1(\kkk; \lll)$. It is easy to see that
\begin{gather}
G_1(\kkk; \lll) \simeq G_1(\lll; \kkk) \simeq G_1(-\kkk; -\lll) \simeq G_1(-\lll; -\kkk) \label{G1_isom_eq} 
\end{gather}
and all other groups of the form $G_1(\kkk;\lll)$ are not isomorphic. 
For a group $\GG_m^2 \rightthreetimes \GG_a$ we use the notation $G_2(\kkk, \lll)$ if it is defined by the character $(t,s) \mapsto t^{-\kkk} s^{-\lll}$ where $(t,s) \in \GG_m^2$ and $k,l \in \ZZ$.

\medskip

\begin{proof}[Proof of~\textup{(\ref{A3_theor_a})}]
Consider a monoid structure $\mu$ on~$\AA^3$, let $X_0$ be the set of invertible elements of~$(\AA^3,\mu)$, and put $G = G(\AA^3,\mu)$ for short. 
Then $G$ is a connected algebraic group with $\dim G = 3$. 
Fix a Levi decomposition $G = L \rightthreetimes U$ where $L$ is a Levi subgroup and $U$ is the unipotent radical of~$G$. 
If $L$ is not a torus, then $\dim L \ge 3$, whence $U$ is automatically trivial and $G = L$ is isomorphic to either $\SL_2$ or $\PSL_2$.
Since $\SL_2$ and $\PSL_2$ have no nontrivial characters, we obtain $X_0 = \AA^3$ by a theorem of Waterhouse~\cite{Wa1982}. 
However, $\SL_2$ and $\PSL_2$ are not isomorphic to $\AA^3$ by~\cite[Lemma~2.1]{KP85}, so this case is excluded. 

In what follows we assume that $L=T$ is a torus, so that $G = T \rightthreetimes U$, and consider several cases depending on the value $\rk(\AA^3, \mu) = \dim T$.

\smallskip

1) $\rk(\AA^3,\mu) = 0$.
As was discussed in Section~\ref{prelim_sec}, in this case $X_0 = \AA^3$ and $(\AA^3,\mu)$ is a unipotent group. It is well known that there exist only two non-isomorphic nilpotent 3-dimensional Lie algebras, so up to isomorphism there exist exactly two 3-dimensional unipotent groups; their multiplications are~$\mathrm{3A}$ and~$\mathrm{U_3}$. 

\smallskip

2) $\rk(\AA^3,\mu) = 3$.
In this case $G$ is a torus. We know from Section~\ref{prelim_sec} that $\mu$ is uniquely determined by the toric structure on~$\AA^3$ and the latter is unique up to isomorphism, so we obtain a unique multiplication on~$\AA^3$, which is~$\mathrm{3M}$. 

\smallskip

3) $\rk(\AA^3,\mu) = 1$.
Since any unipotent group of dimension~$2$ is commutative, we have $U \simeq \GG_a^2$, and so $G \simeq G_1(\kkk;\lll)$ for some $\kkk,\lll \in \ZZ$. Fix a group isomorphism $G \simeq X_0$ and regard $T,U$ as subvarieties of~$X_0$ via this isomorphism. Fix coordinates $t,\alpha,\beta$ in~$X_0$ such that $t\in \GG_m$ is a coordinate in~$T$ and  $\alpha,\beta \in \GG_a$ are coordinates in $U$ satisfying $(t,0,0)*(1,\alpha,\beta)*(t^{-1},0,0) = (1,t^{-\kkk}\alpha ,t^{-\lll}\beta)$.
Then the multiplication $\mu \colon X_0 \times X_0 \to X_0$ is given by
\begin{equation} \label{mult_mu}
(t_1, \alpha_1, \beta_1)*(t_2, \alpha_2,\beta_2) = (t_1t_2, \alpha_2 + t_2^\kkk\alpha_1, \beta_2 + t_2^\lll\beta_1).
\end{equation}
Consider the commutative group $K = T \times U$ and the left commutative reduction $(\AA^3, {}^c\!\mu)$ of $(\AA^3, \mu)$ with respect to~$T$ and~$U$. By Proposition~\ref{comm_red_prop}, the set of invertible elements of $(\AA^3, {}^c\!\mu)$ is $X_0$ and it is isomorphic to~$K$ as an algebraic group. Moreover, identifying $K$ with $X_0$ via the map $(t,(\alpha,\beta)) \mapsto (t,\alpha,\beta)$ we obtain a group embedding $K \hookrightarrow \AA^3$ via the chain $K \xrightarrow{\sim} X_0 \hookrightarrow \AA^3$.

Up to isomorphism, all group embeddings $\iota \colon K \hookrightarrow \AA^3$ were described in the proof of~\cite[Theorem~1]{ABZ2020} and are given by
\begin{equation} \label{emb_of_K}
(t, \alpha, \beta) \mapsto (t, t^b\alpha, t^c(\beta + \kappa\alpha^{d+1})),
\end{equation}
where $b, c, d \in \ZZ$, $\kappa \in \lbrace 0,1 \rbrace$, $\kappa \le b \le c$, and $bd \le c < b(d+1)$. In particular, 
\[b,c \ge 0, \;\; d > 0, \;\; c - bd \ge 0, \ \text{ and } \ c - b(d+1) < 0.\] 
Regarding $K$ as a subset of~$\AA^3$ via~\eqref{emb_of_K}, we need to determine for which values of $b,c,d,\kappa$ the multiplication $\mu$ on~$K$ given by~\eqref{mult_mu} extends to the whole~$\AA^3$. To this end, for $i=1,2,3$ take a point $P_i \in K$ with coordinates $(t_i,\alpha_i,\beta_i)$ in~$K$ and coordinates $(x_i,y_i,z_i)$ in~$\AA^3$ and assume that $P_1*P_2 = P_3$. It remains to express $x_3,y_3,z_3$ via $x_1,y_1,z_1,x_2,y_2,z_2$ and determine the values of $b,c,d,\kappa$ for which these expressions are polynomials.

According to~\eqref{mult_mu}, we have $t_3 = t_1t_2$, $\alpha_3 = \alpha_2 + t_2^k\alpha_1$, $\beta_3 = \beta_2 + t_2^l\beta_1$. Further, in view of~\eqref{emb_of_K}, for each $i = 1,2,3$ there are the relations $x_i = t_i$, $y_i = t_i^b\alpha_i$, $z_i = t_i^c(\beta_i + \kappa \alpha_i^{d+1})$ and $t_i = x_i$, $\alpha_i = x_i^{-b}y_i$, $\beta_i = x_i^{-c}z_i-\kappa x_i^{-b(d+1)}y_i^{d+1}$.
Put also $b' = b+\kkk$, $c' = c+\lll$ for convenience. We compute
\[
x_3 = t_3 = t_1t_2 = x_1x_2
\]
and see that this is always a polynomial. Further,
\[
y_3 = t_3^b\alpha_3 = t_1^bt_2^b(\alpha_2+t_2^k\alpha_1) = x_1^b x_2^b (x_2^{-b}y_2+x_2^kx_1^{-b}y_1) = x_1^by_2 + x_2^{b'}y_1;
\]
the last expression is a polynomial if and only if $b' \ge 0$. At last,
\begin{multline*}
z_3 = t_3^c(\beta_3 + \kappa \alpha_3^{d+1}) = t_1^ct_2^c(\beta_2+t_2^l\beta_1 + \kappa (\alpha_2+t_2^k\alpha_1)^{d+1}) =\\
x_1^c x_2^c (x_2^{-c}z_2-\kappa x_2^{-b(d+1)}y_2^{d+1}+x_2^l(x_1^{-c}z_1-\kappa x_1^{-b(d+1)}y_1^{d+1}) + \kappa (x_2^{-b}y_2+x_2^k x_1^{-b}y_1)^{d+1}) =\\
x_1^cz_2 +x_2^{c'}z_1 + \kappa x_1^c x_2^c((x_2^{-b}y_2+x_2^k x_1^{-b}y_1)^{d+1}-x_2^{-b(d+1)}y_2^{d+1} - x_1^{-b(d+1)}x_2^ly_1^{d+1}).
\end{multline*}
If $\kappa = 0$, then the expression for $z_3$ is a polynomial if and only if $c' \ge 0$. If $\kappa = 1$, then the expression for $z_3$ is a polynomial if and only if $c' \ge 0$ and
\[
Q(x_1,x_2,y_1,y_2) = x_1^c x_2^c((x_2^{-b}y_2+x_2^k x_1^{-b}y_1)^{d+1}-x_2^{-b(d+1)}y_2^{d+1} - x_1^{-b(d+1)}x_2^ly_1^{d+1})
\]
is a polynomial. An easy calculation shows that
\[
Q(x_1,x_2,y_1,y_2) = 
- x_1^{c-b(d+1)}x_2^{c'}y_1^{d+1} + \sum\limits_{\substack{i + j = d+1 \\ 0 \le i \le d}}\!\! \binom{d+1}{i}x_1^{c-bj}x_2^{c-bi+kj}y_1^jy_2^i.
\]
As $c-b(d+1) < 0$, the summand $- x_1^{c-b(d+1)}x_2^{c'}y_1^{d+1}$ must cancel with another term. The latter is possible if and only if this summand cancels with the term of the big sum for $i = 0$ and $j=d+1$, which happens if and only if $\lll = \kkk(d+1)$. Under this condition, all the other terms of the big sum have non-negative exponents if and only if $c+\kkk(d+1)-(b+\kkk)i = c'-b'i \ge 0$ for any $1 \le i \le d$, which holds if and only if $c'-b'd \ge 0$.

Summarizing, we see that for $\kappa = 0$ the expressions of $x_3,y_3,z_3$ are polynomials if and only if $b',c' \ge0$, in which case we get a monoid structure on~$\AA^3$ of type $\mathrm{1M2A}_0(p)$ with $p=(b,b',c,c')$. For $\kappa = 1$ the expressions of $x_3,y_3,z_3$ are polynomials if and only if $b',c' \ge0$, $c'-c = (b'-b)(d+1)$, and $c'- b'd \ge 0$. In particular, we obtain all relations~\eqref{bcbcd_cond_eq}, which imply $b' \ne 0$ and $c' \ne 0$. Thus $p = (b,b',c,c') \in \NN^4$ is a compatible quadruple, $Q(x_1,x_2,y_1,y_2) = Q_p(x_1,x_2,y_1,y_2)$, and we get a monoid structure on~$\AA^3$ of type~$\mathrm{1M2A}_1(p)$.

\smallskip

4) $\rk(\AA^3,\mu) = 2$.
In this case we could refer at once to results of~\cite{Za2024}, but we provide a direct argument using the reduction to the commutative case. The proof is similar to that given in 3) with the following changes: $G \simeq G_2(\kkk,\lll)$ for some $\kkk,\lll \in \ZZ$; we fix coordinates $t,s,\alpha$ in~$X_0$ such that $t,s\in \GG_m$ are coordinates in~$T$ and  $\alpha \in \GG_a$ is a coordinate in~$U$ satisfying $(t,s,0)*(1,1,\alpha) * (t^{-1},s^{-1},0) = (1,1,t^{-\kkk}s^{-\lll}\alpha)$; the multiplication $\mu \colon X_0 \times X_0 \to X_0$ is given by $(t_1,s_1,\alpha_1) * (t_2,s_2,\alpha_2) = (t_1t_2,s_1s_2, \alpha_2 + t_2^{\kkk}s_2^{\lll}\alpha_1)$; the group $K = T \times U$ is identified with $X_0$ via the map $((t,s),\alpha) \mapsto (t,s, \alpha)$, which yields a group embedding $K \hookrightarrow \AA^3$.

Up to isomorphism, all group embeddings $\iota \colon K \hookrightarrow \AA^3$ were described in the proof of~\cite[Proposition~1]{ABZ2020} and have the form $(t,s,\alpha) \mapsto (t,s,t^bs^c\alpha)$ for some $b,c \in \Zgezero$. For $i=1,2,3$ and points $P_i \in K$ with coordinates $(t_i,s_i,\alpha_i)$ in~$K$ and coordinates $(x_i,y_i,z_i)$ in~$\AA^3$ the condition $P_1*P_2 = P_3$ is equivalent to $x_3 = x_1x_2$, $y_3 = y_1y_2$, $z_3 = x_1^by_1^c z_2 + x_2^{b'} y_2^{c'} z_1$ where $b' = b+\kkk$ and $c' = c+\lll$. The latter expression is a polynomial if and only if $b'\ge0$ and $c'\ge0$, in which case we get a monoid structure on~$\AA^3$ of the form~$\mathrm{2M1A}(p)$.
\end{proof}

\begin{proof}[Proof of~\textup{(\ref{A3_theor_b})}]
If two monoids are isomorphic, then their groups of invertible elements are isomorphic as well. So we need to study only monoids of ranks~$1$ and~$2$, i.e., of types $\mathrm{1M2A}_0(p)$, $\mathrm{1M2A}_1(p)$, and $\mathrm{2M1A}(p)$.
Note that for each of these monoids the left and right commutative reductions are computed in Examples~\ref{MMA_red} and~\ref{MAA_red} and they are unique up to isomorphism by Proposition~\ref{choiceTU_prop}.
If two monoids belong to the above list and are isomorphic, then, according to Proposition~\ref{muccmu_prop}, their left commutative reductions are isomorphic and so are their right commutative reductions. Comparing Examples~\ref{MMA_red} and~\ref{MAA_red} with Theorem~\ref{A3comm_theor}(\ref{A3comm_theor_b}) we find that the latter is possible only if the two monoids are of the same type.
So it remains to study which monoids inside each of the three types are isomorphic to each other. 

\smallskip

1) For $p = (b,b',c,c') \in \Zgezero^4$, let us find monoid structures on~$\AA^3$ that are isomorphic to $\mathrm{1M2A}_0(p)$ and of the same type with quadruple $\widetilde p \in \Zgezero^4$. According to Example~\ref{MAA_red}, the left and right commutative reductions of $\mathrm{1M2A}_0(p)$ are isomorphic to $\mathrm{1M2A}_0(b,b,c,c)$ and $\mathrm{1M2A}_0(b',b',c',c')$, respectively.
By Theorem~\ref{A3comm_theor}(\ref{A3comm_theor_b}), a commutative monoid of type $\mathrm{1M2A}_0(b,b,c,c)$ is isomorphic only to itself and $\mathrm{1M2A}_0(c,c,b,b)$. 
So we obtain the following list of possible quadruples~$\widetilde p$:
\begin{enumerate}[label=\textup{(\roman*)},ref=\textup{(\roman*)}]
  \item \label{equ_it} $(b, b', c, c')$;
  \item \label{no1_it} $(b, c', c, b')$;
  \item \label{no2_it} $(c, b', b, c')$;
  \item \label{sym_it} $(c, c', b, b')$.
\end{enumerate}

In cases~\ref{equ_it} and~\ref{sym_it} the quadruples $p, \widetilde p$ are equal and symmetric, respectively, and the respective monoid structures are indeed isomorphic. (In case~\ref{sym_it} the isomorphism is given by $(x,y,z) \mapsto (x,z,y)$.) Let us show that cases~\ref{no1_it} and~\ref{no2_it} are possible only for quadruples $\widetilde p$ that belong to cases~\ref{equ_it} or~\ref{sym_it} as well. Recall from the proof of part~(\ref{A3_theor_a}) that the group of invertible elements of the monoid $\mathrm{1M2A}_0(p)$ is isomorphic to $G_1(\kkk;\lll)$, where $\kkk = b'-b$ and $\lll = c'- c$. Since it is isomorphic to the group of invertible elements of the second monoid with quadruple~$\widetilde p$, in view of~\eqref{G1_isom_eq} the corresponding values $\widetilde \kkk$ and $\widetilde \lll$ satisfy $\{|\kkk|, |\lll|\} = \{|\widetilde\kkk|, |\widetilde\lll|\}$.
Notice that for both cases~\ref{no1_it} and~\ref{no2_it} this means that
\[\{|b'-b|, |c'-c|\} = \{|c'-b|, |b'-c|\}.\]
Then the points $b, b', c, c' \in \ZZ$ are consequent vertices of a (degenerate) kite, whence $b=c$ or $b'=c'$. Swapping these equal variables, we obtain cases~\ref{equ_it} or~\ref{sym_it}.

\smallskip

2) For monoids of type $\mathrm{1M2A}_1(p)$, the beginning of the argument is similar to that in~1) and we again obtain the list \ref{equ_it}--\ref{sym_it} of possible quadruples~$\widetilde p$. Since the quadruple $p$ is compatible, we have $c\ge b$ and $c'\ge b'$. Since $\widetilde p$ is also compatible, all possible quadruples among~\ref{equ_it}--\ref{sym_it} are equal to that in~\ref{equ_it}.

\smallskip

3) Let us determine which monoids of type $\mathrm{2M1A}(p)$, where $p = (b,b',c,c') \in \Zgezero^4$, are isomorphic.
All such monoids have the same set of invertible elements
\[
X_0 = \{x,y \ne 0\} \subseteq \AA^3
\]
whose complement is the union of two irreducible components $C_1 = \{x = 0\}$ and $C_2 = \{y=0\}$. Let us show that the subvarieties
\begin{equation*}
G_{1,1} = \{P_1 \in X_0 \mid P_1 * P_2 = P_2 \text{ for any } P_2 \in C_1\}
\end{equation*}
are not isomorphic for monoids with different values of~$b$. Indeed, for $P_1=(x_1,y_1,z_1) \in \AA^3$ and $P_2=(0,y_2,z_2) \in C_1$ the condition $P_1 * P_2 = P_2$ is equivalent to
\[\begin{aligned}
(y_1-1)y_2 &= 0, \quad (x_1^by_1^c - 1)z_2=0 \quad\text{ if } b' \ne 0;
\\
(y_1-1)y_2 &= 0, \quad (x_1^by_1^c - 1)z_2+y_2^{c'}z_1=0 \quad\text{ if } b' = 0. 
\end{aligned}\]
Thus $P_1 \in G_{1,1}$ if and only if $x_1,y_1 \ne 0$ and
\[\begin{aligned}
y_1 &= 1, \quad x_1^b = 1 \quad\text{ if } b' \ne 0;
\\
y_1 &= 1, \quad x_1^b = 1, \quad z_1=0 \quad\text{ if } b' = 0. 
\end{aligned}\]
Put $\Phi_b = \sqrt[b]{1} = \{\alpha \in \KK \mid \alpha^b = 1\}$; then the variety $G_{1,1}$ is isomorphic to
\[\begin{array}{llll}
\Phi_b \times \KK &\quad\text{ if } b, b' \ne 0; \quad\quad & \KK^\times \times \KK &\quad\text{ if } b = 0,\, b' \ne 0;
\\
\Phi_b &\quad\text{ if } b \ne 0, \, b' = 0; \quad\quad & \KK^\times &\quad\text{ if } b = b' = 0.
\end{array}\]
These varieties are not isomorphic for different values of~$b$. 

In the same way, the subvarieties
\begin{gather*}
G_{1,2} = \{P_1 \in X_0 \mid P_1 * P_2 = P_2 \text{ for any } P_2 \in C_2\},\\
G_{2,1} = \{P_2 \in X_0 \mid P_1 * P_2 = P_1 \text{ for any } P_1 \in C_1\},\\
G_{2,2} = \{P_2 \in X_0 \mid P_1 * P_2 = P_1 \text{ for any } P_1 \in C_2\}
\end{gather*}
are not isomorphic for different values of~$c,b',c'$, respectively. Since any isomorphism between two monoids in question can either preserve $C_1$ and $C_2$ or permute them, the collection $\lbrace (G_{1,1}, G_{2,1}), (G_{1,2}, G_{2,2}) \rbrace$ and hence the collection $\lbrace (b,b'),(c,c') \rbrace$ is an invariant of such monoids under isomorphisms.
It follows that a monoid of type $\mathrm{2M1A}(p)$ with quadruple $p = (b,b',c,c')$ is not isomorphic to another monoid of this type with quadruple~$\widetilde p$ unless $p$ and $\widetilde p$ are equal or symmetric. It remains to notice that in the symmetric case the two monoids are isomorphic via the map $(x,y,z) \mapsto (y,x,z)$.
\end{proof}

The proof of Theorem~\ref{A3_theor} is completed.

\section{Algebraic properties of monoid structures on \texorpdfstring{$\AA^3$}{A\^3}}
\label{properties_sec}

In this section we discuss various algebraic properties of all monoids $(\AA^3,\mu)$, which are classified in Theorem~\ref{A3_theor}. This includes complete descriptions of the groups of invertible elements, centers, idempotents, kernels, and zero elements. We also determine which of the classified monoids are multiplicative monoids of $3$-dimensional associative $\KK$-algebras.

\medskip

\textbf{Groups of invertible elements}.
For each monoid $(\AA^3,\mu)$, Table~\ref{TAB_groups} indicates its set of invertible elements $X_0 \subseteq \AA^3$ and its group of invertible elements $G(\AA^3,\mu)$; see the columns ``$X_0$'' and ``$G(\AA^3,\mu)$'', respectively. (Recall the notations $G_1(\kkk;\lll)$ and $G_2(\kkk;\lll)$ from the beginning of Section~\ref{proof_sec}.) The information is mostly extracted from Examples~\ref{MMA_red} and~\ref{MAA_red} and from the proof of Theorem~\ref{A3_theor}(\ref{A3_theor_a}). Recall that the notation $\mathrm{U_3}$ stands for the group of upper unitriangluar $3\times3$-matrices over~$\KK$.

\begin{table}[h]
 \centering
   \caption{Groups of invertible elements for monoids $(\AA^3,\mu)$}
\label{TAB_groups}
\(\begin{array}{|c|c|c|c|}\hline
\rk & \text{Monoid} & X_0 & G(\AA^3,\mu)
\\\hline\hline
0 & \mathrm{3A} & \AA^3 & \GG_a^3
\\\hline
0 & \mathrm{U_3} & \AA^3 & \mathrm{U_3}
\\\hline
1 & \mathrm{1M2A}_0(b,b',c,c') & \{x \ne 0\} & G_1(b'-b; c'-c)
\\\hline
1 & \mathrm{1M2A}_1(b,b',c,c') & \{x \ne 0\} & G_1(b'-b; c'-c)
\\\hline
2 & \mathrm{2M1A}(b,b',c,c') & \{x,y \ne 0\} & G_2(b'-b, c'-c)
\\\hline
3 & \mathrm{3M} & \{x,y,z \ne 0\} & \GG_m^3
\\\hline
\end{array}\)
\end{table}

\medskip


\textbf{Multiplicative monoids of $3$-dimensional associative algebras}.
Given a finite-dimensional associative $\KK$-algebra~$A$ with unity and of dimension~$k$, the multiplication in~$A$ naturally induces a monoid structure on~$\AA^k$. In the case $k=3$ all such algebras are known from~\cite[Section~4]{Stu1890} and \cite[Section~3]{Sch1891} (see also~\cite[Section~4]{KSTT21}). Up to isomorphism, there are exactly five such algebras; four of them are commutative. In Table~\ref{Mult_monoids_A3} we list all these algebras and provide the corresponding multiplicative monoids, which are easily calculated in each case. The notation $\mathrm{T_2}$ stands for the algebra of upper-triangular $2\times2$-matrices over~$\KK$.

\begin{table}[h]
 \centering
   \caption{Multiplicative monoids of $3$-dimensional algebras}
\label{Mult_monoids_A3}
\(\begin{array}{|c|c|}\hline
\text{Algebra} & \text{Monoid} \\
\hline
\KK\times \KK\times \KK & \mathrm{3M} \\
\hline
\KK\times \KK[t]/(t^2) & \mathrm{2M1A}(0,0,1,1) \\
\hline
\KK[t]/(t^3) & \mathrm{1M2A}_1(1,1,1,1) \\
\hline
\KK[t,s]/(t^2, ts, s^2) & \mathrm{1M2A}_0(1,1,1,1) \\
\hline
\mathrm{T_2} & \mathrm{2M1A}(1,0,0,1) \\
\hline
\end{array}\)
\end{table}

As the first four items in~Table~\ref{Mult_monoids_A3} are the multiplicative monoids of certain commutative (finite-dimensional, associative, and with unity) algebras, we also mention some easy general observations on such monoids.
First, every such algebra uniquely (up to permutation) decomposes into a direct product of local algebras, and so its multiplicative monoid decomposes accordingly.
Second, every local algebra $A$ with maximal ideal~$\mathfrak m$ decomposes (as a variety) as $A = \KK \times \mathfrak m$ with multiplication given by $(x,y)(x',y') = (xx', xy'+x'y+yy')$, and the group of invertible elements of $A$ is $\lbrace (x,y) \mid x \ne 0 \rbrace \simeq \GG_m \times (1 + \mathfrak m)$, where $1 + \mathfrak m$ is commutative and unipotent.

\medskip


\textbf{Centers}.
Recall that the center of a monoid $(X, \mu)$ is the set
\[
Z(X,\mu) = \{x \in X \mid x * y = y * x \ \text{for all} \ y \in X\}.
\]
Clearly, the closure of the center $Z(G(X,\mu))$ of the group $G(X,\mu)$ is a subset of $Z(X,\mu)$. More precisely, the closure of any irreducible component of $Z(G(X,\mu))$ is an irreducible component of~$Z(X,\mu)$. In general, $Z(X,\mu)$ may contain other irreducible components; we call them \emph{external}.

The results of computations of the centers for all monoids $(\AA^3,\mu)$ are gathered in Table~\ref{TAB_center}, where we use the notations
\begin{gather*}
\Phi_m = \sqrt[m]{1} = \{\alpha \in \KK \mid \alpha^m = 1\} 
\quad\quad \text{and} \quad\quad
\Psi_m = \Phi_m \cup \{0\}
\end{gather*}
 for $m \in \NN$ and put $\Phi_0 = \Psi_0 = \KK$ for convenience.

\begin{table}[h]
 \centering
   \caption{The centers of monoids $(\AA^3,\mu)$}
\label{TAB_center}
\(\begin{array}{|c|c|c|c|c|c|}\hline
\rk & \text{Monoid} & Z(\AA^3,\mu) & \dim Z(\AA^3,\mu) & \das\begin{matrix}\text{External}\\\text{components}\end{matrix}
\\\hline\hline
0 & \mathrm{3A} & \AA^3 & 3 & -
\\\hline
0 & \mathrm{U_3} & x=0,\, y=0, \, z \in \KK & 1 & -
\\\hline
\multirow{5}{*}{$1$} & \multirow{5}{*}{$\mathrm{1M2A}_0(b,b',c,c')$}
    & 
    \das\begin{matrix}x \in \Psi_{\gcd(|b'-b|,|c'-c|)} \\ \text{ if } b,b',c,c' \ge 1,\end{matrix}
    & \das\begin{matrix}3 \\ \text{if } b = b', c = c'\end{matrix}
    & -
    \\\cline{4-5}
& & \das \begin{matrix}x \in \Phi_{\gcd(|b'-b|,|c'-c|)} \\ \text{ otherwise};\end{matrix}
    & \das\begin{matrix}1 \\ \text{if } b=b', c \ne c'\end{matrix}
    & \das\begin{matrix}\{(0,y,0)\mid y \in \KK\} \\ \text{ if } b,b',c,c' \ge 1\end{matrix}
    \\\cline{4-5}
& & \das \begin{matrix}y \in \KK \text{ if } b=b', \\ y = 0 \;\text{ otherwise};\end{matrix}
    & \das\begin{matrix}1 \\ \text{if } b \ne b', c = c'\end{matrix}
    & \das\begin{matrix}\{(0,0,z)\mid z \in \KK\} \\ \text{ if } b,b',c,c' \ge 1\end{matrix}
    \\\cline{4-5}
& & \das \begin{matrix}z \in \KK \text{ if } c=c', \\ z = 0 \;\text{ otherwise}\end{matrix}
    & \das\begin{matrix}0 \\ \text{if } b\ne b', c\ne c'\end{matrix}
    & \das\begin{matrix}\{(0,0,0)\} \\ \text{ if } b,b',c,c' \ge 1\end{matrix}
\\\hline
\multirow{2}{*}{$1$} & \multirow{2}{*}{$\mathrm{1M2A}_1(b,b',c,c')$} & \AA^3 \;\text{ if } b=b' & 3 & -
    \\\cline{3-5}
& & \das\begin{matrix}x \in \Psi_{|b'-b|},\, y = 0, \, z = 0 \\\text{ if } b\ne b'\end{matrix} & 0 & \{(0,0,0)\} 
\\\hline
\multirow{4}{*}{$2$} & \multirow{4}{*}{$\mathrm{2M1A}(b,b',c,c')$} & \AA^3 \;\text{ if } (b,c)=(b',c') & 3 & -
    \\\cline{3-5}
& & \das \begin{matrix}x^by^c=x^{b'}y^{c'}, \, z = 0 \\\text{ if } (b,c) \ne (b',c')\end{matrix} 
    & 1 & \das\begin{matrix}
\{(0,y,0)\mid y \in \KK\} \\ \text{ if }b,b' \ge 1,\\
\{(x,0,0)\mid x \in \KK\} \\ \text{ if }c,c' \ge 1
\end{matrix}
\\\hline
3 & \mathrm{3M} & \AA^3 & 3 & -
\\\hline
\end{array}\)
\end{table}

An inspection of Table~\ref{TAB_center} shows that the center $Z(\AA^3,\mu)$ is always an unmixed variety, i.e., has equidimensional irreducible components. We do not know whether this holds for all monoid structures on~$\AA^n$ with arbitrary~$n$.
It is proved in~\cite[Corollary~3]{Za2024} that for noncommutative monoid structures of rank $n-1$ on $\AA^n$ the dimensions of all irreducible components of $Z(\AA^n,\mu)$ are equal to~$n-2$. Further, it follows from Table~\ref{TAB_center} that for noncommutative monoid structures on~$\AA^3$ these dimensions are always at most~1.

At the same time, there is an example of a monoid structure~$\mu$ of rank~$2$ on the variety $X = \{vw=zt\} \subseteq \AA^4$ of dimension~$3$ such that $Z(X,\mu)$ has irreducible components of dimensions $0$ and~$1$; see \cite[Example~7]{Za2024}. 

As can be seen from Table~\ref{TAB_center}, certain noncommutative monoids $(\AA^3,\mu)$ of rank~1 have a finite center. Notice that the order of such center can be equal to any positive integer. 

Finally, let us consider ``the most noncommutative'' monoid structures, i.e., those with trivial (one-element) center. On $\AA^2$ there is only one such structure up to isomorphism; it is given by
\[(x_1, y_1) * (x_2, y_2) = (x_1x_2, y_1+x_1y_2).\]
The corresponding monoid is isomorphic to the monoid $\left\{\das\begin{pmatrix}x & y \\ 0 & 1\end{pmatrix} \mid x,y \in \KK\right\}$ with matrix multiplication. 
At the same time, there are infinitely many pairwise non-isomorphic monoid structures on~$\AA^3$ with trivial center: this is the case of type $\mathrm{1M2A}_0(b,b',c,c')$, where $b \ne b'$, $c \ne c'$, at least one of $b,b',c,c'$ equals~$0$, and $\gcd(|b'-b|, |c'-c|) = 1$. For example, for $b=c=1$ and $b'=c'=0$ we obtain the monoid structure given by
\[(x_1, y_1, z_1) * (x_2, y_2, z_2) = (x_1x_2, y_1+x_1y_2, z_1 +x_1z_2),\]
which corresponds to the monoid $\left\{\das\begin{pmatrix} x & y & z \\ 0 & 1 & 0 \\ 0 & 0 & 1 \end{pmatrix} \mid x,y,z \in \KK\right\}$ with matrix multiplication.

\medskip


\textbf{Idempotents}.
Given a monoid~$(X, \mu)$, recall that an element $e \in X$ is called an \emph{idempotent} if $e * e = e$. The set $E(X,\mu)$ of all idempotents is clearly a closed subvariety in~$X$. 
Let
\[E(X,\mu) = \bigcup\limits_i E_i\]
be the decomposition of $E(X,\mu)$ into the union of irreducible components. We assume that the component $E_1 := E(X,\mu) \cap G(X,\mu) = \{{\bf 1}\}$ is the first one. 

For any $a, b \in \ZZ$, denote by $z(a,b)$ the number of zeros among $a$ and~$b$, i.e., 
\[z(a,b) = \begin{cases}
2 \quad \text{ if } a=0, b=0; \\
0 \quad \text{ if } a \ne 0, b \ne 0;\\
1 \quad \text{ otherwise.}
\end{cases}\]

By direct computations, we find the set of idempotents for each of the monoids $(\AA^3,\mu)$. The results are presented in Table~\ref{TAB_idemp}. 

\begin{table}[h]
 \centering
   \caption{Irreducible components of the subvariety $E(\AA^3,\mu)$}
\label{TAB_idemp}
\(\begin{array}{|c|c|c|c|c|c|}\hline
\rk & \text{Monoid} & i & \text{Case} & E_i & \dim E_i
\\\hline\hline
0 & \mathrm{3A} & 1 & & \{(0,0,0)\} & 0 
\\\hline
0 & \mathrm{U_3} & 1 & & \{(0,0,0)\} & 0 
\\\hline
1 & \mathrm{1M2A}_0(b,b',c,c') & 1 & & \{(1,0,0)\} & 0
    \\\cline{3-6}
& & 2 & z(b,b')=1, \, z(c,c')=1 & \{(0,y,z)\mid y,z \in \KK\} & 2
    \\\cline{4-6}
& & & z(b,b')\ne 1, \, z(c,c')=1 & \{(0,0,z)\mid z \in \KK\} & 1
    \\\cline{4-6}
& & & z(b,b')=1, \, z(c,c')\ne1 & \{(0,y,0)\mid y \in \KK\} & 1
    \\\cline{4-6}
& & & z(b,b')\ne1, \, z(c,c')\ne1 & \{(0,0,0)\} & 0
\\\hline
1 & \mathrm{1M2A}_1(b,b',c,c') & 1 & & \{(1,0,0)\} & 0
    \\\cline{3-6}
& & 2 & & \{(0,0,0)\} & 0
\\\hline
2 & \mathrm{2M1A}(b,b',c,c') & 1 & & \{(1,1,0)\} & 0
    \\\cline{3-6}
& & 2 & z(b,b')=1 & \{(0,1,z)\mid z \in \KK\} & 1
    \\\cline{4-6}
& & & \text{otherwise} & \{(0,1,0)\} & 0
    \\\cline{3-6}
& & 3 & z(c,c')=1 & \{(1,0,z)\mid z \in \KK\} & 1
    \\\cline{4-6}
& & & \text{otherwise} & \{(1,0,0)\} & 0
    \\\cline{3-6}
& & 4 & \das
    \begin{matrix}\text{exactly one of }(b,c),\\(b',c') \text{ equals }(0,0)\end{matrix} 
    & \{(0,0,z)\mid z \in \KK\} & 1
    \\\cline{4-6}
& & & \text{otherwise} & \{(0,0,0)\} & 0
\\\hline
3 & \mathrm{3M} & \text{1--8} & & \{(x, y, z) \mid x,y,z =0,1\} & 0
\\\hline
\end{array}\)
\end{table}

Let us mention some general facts on idempotents in algebraic monoids. If $\bf 1$ is the unique idempotent of a monoid~$(X,\mu)$, then $X=G(X,\mu)$; see~\cite[Proposition~3~(iii)]{Br-1}. This agrees with Table~\ref{TAB_idemp}, where the idempotent is unique only for the cases $\mathrm{3A}$ and~$\mathrm{U_3}$. 

Assume that $X$ is irreducible. By~\cite[Theorem~1.1]{Br-2}, the variety $E(X,\mu)$ is smooth (and hence a disjoint union of its irreducible components), and each irreducible component $E_i$ is an orbit for the action of~$G(X,\mu)$ on~$X$ by conjugation. In particular, the isolated idempotents are exactly the central idempotents, and the irreducible components $E_i$ with $\dim E_i > 0$ do not intersect the center $Z(X,\mu)$. For monoids studied in this paper, one can check this by inspecting Tables~\ref{TAB_center} and~\ref{TAB_idemp}.

From Table~\ref{TAB_idemp} one can see that for every monoid $(\AA^3,\mu)$ the number of irreducible components in $E(\AA^3,\mu)$ equals~$2^r$, where $r$ is the rank of the monoid. So it is natural to pose the following question.

\begin{question}
\label{quest_idemp}
Consider any monoid structure $\mu$ on~$\AA^n$ of rank~$r$. Does the number of irreducible components in $E(\AA^n,\mu)$ equal~$2^r$? 
\end{question}

The answer to Question~\ref{quest_idemp} is affirmative provided $n \le 3$. Indeed, the case $n=3$ follows from Table~\ref{TAB_idemp}. For $n = 2$ and $r=1$ this is true thanks to~\cite[Theorem~3]{Za2024}; see also~\cite[Example~5]{Za2024}. Other cases with $n \le 2$ are either coordinatewise addition or coordinatewise multiplication.

Given a monoid $(X,\mu)$ with irreducible affine~$X$, consider a maximal torus $T$ in the group $G(X,\mu)$, regard~$T$ as a subvariety of~$X$, and let $\overline T$ be the closure of~$T$ in~$X$, which is a (possibly non-normal) affine toric variety. According to~\cite[Theorem 2.14]{Br-2}, every irreducible component of $E(X,\mu)$ meets~$\overline{T}$. Note that $\overline T$ is a submonoid of~$X$; then by~\cite[Theorem IV.9]{Neeb} every $T$-orbit in $\overline T$ contains a unique idempotent. It follows that the number of irreducible components in $E(X,\mu)$ does not exceed the number of $T$-orbits in $\overline{T}$, which equals the number of faces in the cone $C(\overline T)$ assigned to $\overline T$ as an affine toric variety.

\begin{question} \label{quest_T-orbits}
Consider any monoid structure $\mu$ on~$\AA^n$. Is it true that every irreducible component of $E(\AA^n,\mu)$ meets exactly one $T$-orbit in~$\overline T$?
\end{question}

For commutative monoids $(\AA^n,\mu)$, the answer to Question~\ref{quest_T-orbits} is positive. Indeed, in this case all idempotents are central and hence isolated, so every irreducible component of $E(\AA^n,\mu)$ consists of a single point, which belongs to only one $T$-orbit in~$\overline T$. In particular, $|E(\AA^n,\mu)|$ coincides with the number of faces of the cone $C(\overline T)$. It was proved in~\cite[Proposition~10]{ABZ2020} that the latter number is always at least $2^r$ and equals $2^r$ for $r \le 2$, where $r = \rk(\AA^n,\mu)$. In view of the above and Question~\ref{quest_idemp}, the following question seems to be natural. (Notice that part~(b) is a stronger version of part~(a).)

\begin{question} \label{quest_cones}
Let $\mu$ be a monoid structure on~$\AA^n$ of rank~$r$. \\
(a) Does the cone $C(\overline{T})$ have $2^r$ faces?\\
(b) Is the variety $\overline{T}$ isomorphic to $\AA^r$? 
\end{question}

One can extract from the third column of Table~\ref{TAB_groups} that the answer to Question~\ref{quest_cones}(b) (and hence to Question~\ref{quest_cones}(a)) is positive for all monoid structures on~$\AA^3$. Easy calculations show that the same holds for all monoid structures on~$\AA^2$ and~$\AA^1$.

It is known that any monoid $(X,\mu)$ with irreducible affine~$X$ and nilpotent group $G(X,\mu)$ has finitely many idempotents, but this does not extend to the case of solvable $G(X,\mu)$; see~\cite[Theorem~1.12 and Example~1.15]{Pu-81}. Our results also confirm the last observation; see Table~\ref{TAB_idemp}. More results on affine monoids with finitely many idempotents can be found in~\cite{Hu96a, Hu96b}.

\medskip


\setlength{\arraycolsep}{2pt}

\begin{table}[h]
 \centering
   \caption{Kernels and zeros for monoids $(\AA^3,\mu)$}
\label{TAB_zeros}
\(\begin{array}{|c|c|c|c|c|}\hline
\rk & \text{Monoid} & ((0,\!0,\!0)\hspace{0.5pt}{*}\hspace{0.5pt}\AA^3\hspace{0.5pt}{*}\hspace{0.5pt}(0,\!0,\!0),\,\mu) & \text{Kernel} & \text{Zero}
\\\hline\hline
0 & \mathrm{3A} & \mathrm{3A}=(\AA^3,+) & \AA^3 & \text{no}
\\\hline
0 & \mathrm{U_3} & \mathrm{U_3} = (\AA^3,*) & \AA^3 & \text{no}
\\\hline
1 & \mathrm{1M2A}_0(b,b',c,c') & 
  \das\begin{matrix}
  (\{(0,y,z)\}, +) \\ \text{if } b,b',c,c'=0;\,\quad\quad\quad \\
  (\{(0,y,0)\}, +) \\ \text{if } b,b'=0,\;\; c+c'\ge1; \\
  (\{(0,0,z)\}, +) \\ \text{if } b+b'\ge1,\;\; c,c'=0; \\
  \{(0,0,0)\} \\ \text{if } b+b'\ge1,c+c'\ge1
  \end{matrix} & \setlength{\arraycolsep}{2pt} \das\begin{array}{cl}
  \{(0,y,z)\} & \text{if } z(b,b')\ge1 \ \text{and}\\
                  & z(c,c')\ge1; \\
  \{(0,y,0)\} & \text{if } b, b' \ge 1 \ \text{and}\\
                  & z(c,c')\ge1; \\
  \{(0,0,z)\} & \text{if } z(b,b')\ge1 \ \text{and}\\
                  & c,c'\ge1; \\
  \{(0,0,0)\} & \text{if } b,b',c,c'\ge1
  \end{array} & \das\begin{matrix}
  \text{yes iff } \\ b,b', c,c' \ge 1
  \end{matrix}
\\\hline
1 & \mathrm{1M2A}_1(b,b',c,c') & \{(0,0,0)\} & \{(0,0,0)\} & \text{yes}
\\\hline
2 & \mathrm{2M1A}(b,b',c,c') &
  \das\begin{matrix}
  (\{(0,0,z)\}, +) \\ \text{if } b,b',c,c'=0;\\
  \{(0,0,0)\} \\ \text{otherwise}
  \end{matrix} & \das\begin{matrix}
  \{(0,0,z)\} \\ \text{if } (b,c) \text{ or } (b',c') \text{ equals } (0,0);\\
  \{(0,0,0)\} \\ \text{if } (b,c), (b',c')\ne(0,0)
  \end{matrix} & \das\begin{matrix}
  \text{yes iff } \\ (b,c) \ne (0,0), \\ (b',c') \ne (0,0)
  \end{matrix}
\\\hline
3 & \mathrm{3M} & \{(0,0,0)\} & \{(0,0,0)\} & \text{yes}
\\\hline
\end{array}\)
\end{table}

\textbf{Kernels and zero elements}.
Given a monoid $(X,\mu)$, recall that an element ${\bf 0} \in X$ is called a \emph{zero} if ${\bf 0} * x = x * {\bf 0} = {\bf 0}$ for any $x \in X$, and a subset $I \subseteq X$ is called an \emph{ideal} if $X*I*X \subseteq I$. Clearly, a zero element ${\bf 0} \in X$ is unique if exists, and in this case the subset $\lbrace \bf 0 \rbrace$ is automatically an ideal. If $X$ is irreducible, then by~\cite[Theorem~1]{Ri1} there exists the smallest (with respect to inclusion) ideal; it is called the \emph{kernel} of~$X$. It follows that $(X,\mu)$ admits a zero if and only if its kernel consists of one element (which is necessarily a zero). In Table~\ref{TAB_zeros}, we present the kernels of all monoids $(\AA^3,\mu)$ and mark in which cases there exists a zero. 

For computing the kernels, we use the following observations. There is the classical partial order on the set of idempotents: $e_1 \le e_2$ if $e_1*e_2 = e_2*e_1 = e_1$. It is known that an idempotent $e$ is minimal with respect to this order if and only if the set $e*X*e$ is a group with respect to the multiplication~$\mu$; see~\cite[Proposition~4]{Br-1}. In each case, we compute the set $e*X*e$ for the idempotent $e=(0,0,0)\in X=\AA^3$ and find out that it is a group (see the third column of Table~\ref{TAB_zeros}), so $(0,0,0)$ is a minimal idempotent. By~\cite[Proposition~5(v)]{Br-1}, if $e$ is a minimal idempotent in~$X$, then the kernel of~$X$ equals $X*e*X$, so all kernels are found as $\AA^3 * (0,0,0) * \AA^3$. 

The following question has a positive answer for all $n \le 3$.

\begin{question}
Let $\mu$ be a monoid structure on~$\AA^n$ such that $(\AA^n,\mu)$ is not isomorphic to the direct product of two monoids with underlying varieties of positive dimension. Suppose that the center $Z(\AA^n,\mu)$ has an external component. Is it true that $(\AA^n,\mu)$ admits a zero?
\end{question}

\appendix

\section*{Appendix. Reductive monoid structures on affine spaces}

In this appendix, which is independent of the main part of the paper, we provide a classification of all reductive monoid structures on affine spaces of arbitrary dimension together with a short proof; see Theorem~\ref{thm_redmon} below. The result can be also deduced from the classification of smooth reductive monoids given by~\cite[Theorem~27.25]{Ti2011}; see~\cite[Corollary~8.4.5]{Re1985} for a particular case.

For every $n \in \ZZ_{>0}$ let $\mathrm M_n$ denote the set of all $n\times n$-matrices over~$\KK$ and let $\nu_n \colon \mathrm M_n \times \mathrm M_n \to \mathrm M_n$ be the matrix multiplication. Note that $(\mathrm M_n, \nu_n)$ is a reductive monoid with group of invertible elements~$\GL_n$. For $n_1,\ldots,n_k \in \ZZ_{>0}$ put $\mathrm{R}(n_1,\ldots,n_k) = \prod\limits_{i=1}^k (\mathrm M_{n_i},\nu_{n_i})$; this is a reductive monoid with group of invertible elements $\prod\limits_{i=1}^k \GL_{n_i}$.

\begin{theorem} \label{thm_redmon}
The following assertions hold.
\begin{enumerate}[label=\textup{(\alph*)},ref=\textup{\alph*}]
\item \label{thm_redmon_a}
Every reductive monoid of the form $(\AA^n,\mu)$ is isomorphic to $\mathrm R(n_1,\ldots, n_k)$ for some $n_1,\ldots,n_k \in \ZZ_{>0}$.
\item \label{thm_redmon_b}
Two monoids $\mathrm R(n_1,\ldots, n_k)$ and $\mathrm R(m_1,\ldots, m_l)$ are isomorphic if and only if $k = l$ and the tuple $(m_1,\ldots, m_l)$ is obtained from $(n_1,\ldots,n_k)$ by a permutation.
\end{enumerate}
\end{theorem}

In the proof we will use the following notation: for an algebraic group~$G$ and a finite-dimensional $G$-module~$V$, the corresponding dual $G$-module is denoted by~$V^*$.

\begin{proof}[Proof of Theorem~\textup{\ref{thm_redmon}}]
(\ref{thm_redmon_a})
Let $\mu$ be a monoid structure on~$\AA^n$, put $G = G(\AA^n,\mu)$, and suppose that $G$ is reductive.
Consider the action of $G \times G$ on $\AA^n$ given by $(g_1,g_2) \cdot x = g_1 * x * g_2^{-1}$ and the corresponding induced $G\times G$-module structure on~$\KK[\AA^n]$.
Recall from Section~\ref{prelim_sec} that $G \subseteq \AA^n$ is an open $G\times G$-orbit, and so there is an inclusion
\begin{equation} \label{inclusion_eq}
\KK[\AA^n] \subseteq \KK[G]
\end{equation}
of $G \times G$-modules. Further, it is well known (see, for instance, \cite[Theorem~2.15]{Ti2011}) that there is a $G\times G$-module isomorphism
\begin{equation} \label{K[G]_eq}
\KK[G] \simeq \bigoplus \limits_{V} V^*\otimes V,
\end{equation}
where $V$ runs over all pairwise non-isomorphic simple finite-dimensional $G$-modules and the left (resp. right) factor of $G\times G$ acts on the left (resp. right) tensor factor of $V^* \otimes V$.

Since $\AA^n$ has an open $G\times G$-orbit, it follows that all $G\times G$-invariant functions in~$\KK[\AA^n]$ are constants.
Then by~\cite[Proposition~5.1]{KP85} the action of~$G\times G$ on~$\AA^n$ is linearizable.
In what follows we assume that $G\times G$ acts linearly on~$\AA^n$ and fix a decomposition $\AA^n = \bigoplus \limits_{i=1}^k W_i$ into a direct sum of simple $G\times G$-modules. For every $i = 1,\ldots,k$ we have $W_i^* \subseteq \KK[\AA^n]$, therefore by~(\ref{inclusion_eq}) and~(\ref{K[G]_eq}) there is a simple $G$-module $V_i$ such that $W_i^* \simeq V_i^* \otimes V_i$ and hence $W_i \simeq V_i \otimes V_i^*$ as $G\times G$-modules. Then the action of $G$ on~$\AA^n$ by left multiplication is given by its actions on the left tensor factors of each $V_i \otimes V_i^*$ and hence is determined by a homomorphism $\varphi \colon G \to \prod \limits_{i=1}^k \GL(V_i)$. Since this action admits an open orbit (which is again~$G$) with trivial stabilizers and $\dim (\bigoplus \limits_{i=1}^k W_i) = \dim (\prod \limits_{i=1}^k \GL(V_i))$, we conclude that~$\varphi$ is an isomorphism and thus the monoid $(\AA^n,\mu)$ is isomorphic to $\mathrm{R}(n_1, \ldots, n_k)$, where $n_i = \dim V_i$ for all $i=1,\ldots,k$.

(\ref{thm_redmon_b})
If two monoids $\mathrm R(n_1,\ldots, n_k)$ and $\mathrm R(m_1,\ldots, m_l)$ are isomorphic, then so are the groups $\prod\limits_{i=1}^k \GL_{n_i}$ and $\prod\limits_{j=1}^l \GL_{m_j}$, which is possible only if $k = l$ and the tuple $(m_1,\ldots, m_l)$ is obtained from $(n_1,\ldots,n_k)$ by a permutation. Conversely, under the latter conditions the monoids $\mathrm R(n_1,\ldots, n_k)$ and $\mathrm R(m_1,\ldots, m_l)$ are clearly isomorphic.
\end{proof}

\begin{remark}
The Wedderburn--Artin theorem asserts that every semisimple finite-dimensional associative $\KK$-algebra is isomorphic to the direct product $\mathrm M_{n_1} \times \ldots \times \mathrm M_{n_k}$, where each $\mathrm M_{n_i}$ is regarded as the algebra of $n_i\times n_i$-matrices over~$\KK$ (see, for instance,~\cite[Section~3.5, Corollary~b]{Pie82}). It follows from Theorem~\ref{thm_redmon} that, up to isomorphism, all reductive monoids of the form $(\AA^n,\mu)$ are precisely the multiplicative monoids of semisimple finite-dimensional associative $\KK$-algebras.
\end{remark}

\end{document}